\documentclass{opt2026}

\usepackage{amsmath,amssymb}
\usepackage{array}
\usepackage{booktabs}
\usepackage{graphicx}
\usepackage{microtype}
\usepackage{multirow}
\usepackage{placeins}
\usepackage{float}
\usepackage{siunitx}
\usepackage{tikz}
\usetikzlibrary{arrows.meta,calc,positioning}

\newcommand{\R}{\mathbb{R}}
\newcommand{\proj}{\operatorname{proj}}

\newcommand{\nnz}{\operatorname{nnz}} 
\newcommand{\solvername}{\textsc{ShardLP}}

\usepackage{colortbl}
\usepackage[normalem]{ulem}
\definecolor{arxivgreen}{HTML}{008000}
\definecolor{sharddark}{HTML}{0F6F6B}
\definecolor{shardfill}{HTML}{E9F5F4}
\definecolor{cpudark}{HTML}{9A5A13}
\definecolor{cpufill}{HTML}{FFF4E5}

\newif\ifarxivchanges
\arxivchangesfalse
\DeclareRobustCommand{\arxivchange}[1]{%
  \begingroup
  \ifarxivchanges
    \color{blue}%
    \hypersetup{urlcolor=blue,linkcolor=blue,citecolor=blue}%
  \fi
  #1%
  \endgroup
}
\DeclareRobustCommand{\arxivremoved}[1]{%
  \ifarxivchanges
    \leavevmode
    {\color{red}\hypersetup{urlcolor=red,linkcolor=red,citecolor=red}\sout{#1}}%
  \fi
}
\DeclareRobustCommand{\arxivremovedparagraph}[1]{%
  \ifarxivchanges
    \begingroup\sloppy
    \arxivremoved{#1}\par
    \endgroup
  \fi
}
\DeclareRobustCommand{\arxivreplaced}[2]{%
  \ifarxivchanges
    \arxivremoved{#1}\hspace{0.2em}%
  \fi
  \arxivchange{#2}%
}

\DeclareRobustCommand{\arxivgreen}[1]{%
  \begingroup
  \ifarxivchanges
    \color{arxivgreen}%
    \hypersetup{urlcolor=arxivgreen,linkcolor=arxivgreen,citecolor=arxivgreen}%
  \fi
  #1%
  \endgroup
}
\DeclareRobustCommand{\arxivgreenreplace}[2]{%
  \ifarxivchanges
    \arxivremoved{#1}\hspace{0.2em}%
  \fi
  \arxivgreen{#2}%
}

\makeatletter
\renewcommand*{\@titlefoot}{}
\makeatother

\title[Distributed LP on GPU Clusters]{Distributed Linear Programming on GPU
Clusters at Extreme Scale}

\optauthor{%
\Name{Arnaud Deza}\textsuperscript{*}
\Email{adeza3@gatech.edu}\\
\Name{Santanu S. Dey}
\Email{santanu.dey@isye.gatech.edu}\\
\Name{Pascal Van Hentenryck}
\Email{pvh@gatech.edu}\\
\addr H. Milton Stewart School of Industrial and Systems Engineering\\
Georgia Institute of Technology, Atlanta, Georgia, USA
}

\begin{document}
\maketitle

\begingroup
\renewcommand{\thefootnote}{*}
\renewcommand{\footnoteseptext}{\ }
\footnotetext{Corresponding author: \texttt{adeza3@gatech.edu}}
\endgroup

\begin{abstract}
Large linear programs can exceed the memory of a single compute node. Although first-order methods replace sparse factorizations with GPU-suited matrix--vector products, other solver phases can reintroduce a single-node memory limit.  We present \solvername, a distributed GPU LP solver that keeps the matrix and primal--dual state partitioned from sharded input through solution output.  On the Google PDLP benchmark, \solvername{} \arxivreplaced{reaches the published criterion on nine}{solves ten} of eleven instances, compared with eight in the published CPU PDLP study. On the largest benchmark, eight H200 GPUs solve a 1.185-billion-variable, 6.338-billion-nonzero LP in 9.9 minutes; the published CPU experiment reports 21.06 hours on different hardware. Beyond this benchmark, separately checked multi-node solves reach up to 13.604 billion variables and 40.807 billion nonzeros, while validated executions span up to 76 GPUs across 29 compute nodes. For column-partitioned solves, support-aware communication skips GPUs that store no coefficients for a row; on an LP with 2.76 billion nonzeros, it cuts modeled communication by 92.97\% and improves solver time by $1.27\times$--$1.52\times$.
\arxivchange{The code for \solvername{} is available at \url{https://github.com/ArnaudDeza/ShardLP}.}
\end{abstract}

\section{Introduction}

Linear programming (LP) is a core model in large-scale optimization. Modern simplex and interior-point methods are highly effective, but at very large scale their factorization-based linear algebra can become memory intensive and difficult to parallelize \citep{applegate2026pdlp}. PDHG-based first-order methods have a different computational profile: their dominant operations are sparse matrix--vector products, projections, and vector operations, which map naturally to GPUs and distributed memory.

Primal--dual hybrid gradient (PDHG), also known as the Chambolle--Pock method, is the first-order method underlying Primal-Dual Linear Programming (PDLP) \citep{chambolle2011first,applegate2021practical}.  PDLP combines the basic iteration with diagonal scaling, presolve, adaptive step sizes, restart, and feasibility polishing \citep{applegate2023faster,applegate2026pdlp}.  We use the eleven released instances from the large-scale Google PDLP study \citep{applegate2026pdlp} as our principal benchmark. Its largest instance contains 1.185 billion variables and 6.338 billion nonzeros.

Recent GPU implementations show that PDHG maps well to accelerators, but scaling a general sparse LP solver across multiple compute nodes remains much less developed. Among the closest systems, D-PDLP \citep{li2026dpdlp} \arxivreplaced{distributes the two PDHG matrix products across eight H100 GPUs within a single node}{distributes PDHG across eight H100 GPUs within a single node}, while MPAX \citep{lu2024mpax} leaves efficient distributed sparse-data sharding as future work. Our question is therefore not only how to distribute the PDHG iteration, but how to turn it into a complete solver when neither the matrix nor the primal--dual state fits on one node. A broader comparison with prior GPU and distributed optimization systems is given in Appendix~\ref{app:scope}.
 
We make three contributions. First, \solvername{} keeps the matrix and primal--dual state distributed \arxivreplaced{during every solver stage---\mbox{\emph{persistent ownership}}---so no phase needs a full copy on one node}{during scaling, optimization, polishing, recovery, and output---\emph{persistent ownership}---so that the memory needed by these stages is set by the local matrix blocks and vector slices}. Second, it \arxivreplaced{reaches the Google PDLP study's published criterion on nine of eleven benchmark instances}{solves ten of eleven Google benchmark instances} and separately validates multi-node solves up to 13.604 billion variables and 40.807 billion nonzeros. Third, we show that communication can erase multi-GPU speedups and introduce support-aware communication, which skips GPUs with no coefficients for a row. On Design Match, this cuts modeled communication by 92.97\% and improves solver time by $1.27\times$--$1.52\times$.

\arxivchange{The code for \solvername{} is available at \url{https://github.com/ArnaudDeza/ShardLP}.}

\section{From a distributed iteration to a distributed solver}

\subsection{PDHG and the D-PDLP matrix decomposition}

Consider a linear program in interval form, with $c\in\R^n$,
\begin{equation}
  \min_{x\in X} c^\top x
  \qquad \text{s.t.}\qquad Ax\in S,
  \quad
  X=[\ell_v,u_v],\quad S=[\ell_c,u_c],\quad A\in\R^{m\times n}.
  \label{eq:lp}
\end{equation}
Here, $X$ gives the variable bounds and $S$ the lower and upper bounds on each row activity.  At a high level, PDHG alternates a primal update that uses $A^\top y$ and a dual update that uses $Ax$.  \arxivgreenreplace{Suppressing diagonal preconditioning, following \mbox{\citep{chambolle2011first,applegate2021practical}}, define the base map}{Suppressing diagonal preconditioning, we define the PDHG map \citep{chambolle2011first,applegate2021practical}} $T_{\tau,\sigma}(w)=(\widehat x,\widehat y)$ for $w=(x,y)$ by
\begin{align}
 \widehat x &= \proj_X\!\left(x-\tau(c-A^\top y)\right),
 \label{eq:pdhg-primal}\\[-2pt]
 \widehat y &= \Phi_{S,\sigma}\!\left(y,
                 A(2\widehat x-x)\right). \label{eq:pdhg-dual}
\end{align}
Here $\tau$ and $\sigma$ are the primal and dual step sizes, $\proj_X$ clips to the variable bounds, and $\Phi_{S,\sigma}$ denotes the coordinatewise dual update for the row bounds.  \solvername{} does not introduce a new PDHG variant. Its optimization logic follows cuPDLPx \citep{lu2025cupdlpx}; within each restart epoch, it applies the reflected-Halpern update
\begin{equation}
 w^{k+1}=\frac{k+1}{k+2}\bigl((1+\gamma)T_{\tau,\sigma}(w^k)-\gamma w^k\bigr)
          +\frac{1}{k+2}w^0,\qquad \gamma\in[0,1],
 \label{eq:halpern}
\end{equation}
where $\gamma$ controls the reflection and a restart resets the anchor $w^0$ and epoch counter. Scaling, step selection, primal-weight adaptation, and restart follow \citet{lu2025cupdlpx}; their control logic uses scalar reductions.

To distribute the two matrix products, \solvername{} adopts D-PDLP's two-dimensional matrix decomposition \citep{li2026dpdlp}, with one Message Passing Interface (MPI) process, or \emph{rank}, per GPU.  We arrange the $p=RC$ ranks in an $R\times C$ process grid, with $R$ constraint (matrix-row) partitions and $C$ variable (matrix-column) partitions.  For each experiment, we choose the \(R\times C\) grid \arxivreplaced{before sharding}{before partitioning} so that the estimated matrix block and replicated vector slices fit in each GPU's memory.  Row set $I_r$ and column set $J_c$ define the block $A_{rc}=A_{I_r,J_c}$.  The corresponding vector blocks are $y_r=y_{I_r}$ and $x_c=x_{J_c}$.  The two products are
\begin{equation}
 \left(A^\top y\right)_{J_c}=\sum_{r=1}^{R} A_{rc}^\top y_r,
 \qquad
 \left(Ax\right)_{I_r}=\sum_{c=1}^{C} A_{rc}x_c.
 \label{eq:distributed-products}
\end{equation}
Rank $(r,c)$ stores $A_{rc}$.  Ranks sharing $c$ form a \emph{process column}: they replicate $x_c$ and combine partial $A^\top y$ products.  Ranks sharing $r$ form a \emph{process row}: they replicate $y_r$ and combine partial $Ax$ products.  \arxivgreenreplace{Storage per rank is}{For balanced row and column partitions, storage per rank is} $O(\nnz(A_{rc})+n/C+m/R)$.  \arxivreplaced{Up to this point, the decomposition is D-PDLP's. \mbox{\solvername{}}'s contribution is to preserve this ownership beyond the two matrix products and throughout the rest of the solve, as illustrated in Figure~\mbox{\ref{fig:persistent-sharding}}.}{D-PDLP also distributes Halpern updates, convergence checks, restarts, and primal-weight adaptation. \solvername{} retains these computations.}

\begin{figure}[t]
  \centering
  \input{figures/persistent_sharding.tex}
  \caption{Why persistent ownership matters. (a) Centralized solver phases impose an $O(\nnz(A)+m+n)$ per-node memory ceiling even when PDHG is distributed. (b) In a $2\times3$ process grid, rank $(r,c)$ stores $A_{rc}$; row/column reductions form $Ax$ and $A^\top y$ \arxivreplaced{while ownership remains sharded from input through output}{and scaling, polishing, recovery, and output operate on the same local blocks and vector slices}.}
  \label{fig:persistent-sharding}
\end{figure}

\subsection{Persistent ownership beyond the iteration}
\arxivremovedparagraph{Distributing only $Ax$ and $A^\top y$ is not enough: if scaling, polishing, recovery, or another solver stage gathers the full problem on one node, that node still sets the memory limit. \mbox{\solvername{}} therefore keeps the same partition throughout the run. Each rank reads a prepartitioned \emph{prepared shard}; constructing these shards from a monolithic model is an offline step excluded from timing. Scaling, restart, termination, polishing, recovery, and output remain distributed. General presolve is disabled. A distributed pass fixes variables forced to zero by globally singleton zero-equality rows; in both arms of the RQ3 Design Match communication experiments, this removes 1.123 billion nonzeros (40.69\%), leaving 1.637 billion nonzeros. The RQ1 experiments do not use this reduction. Appendix~\mbox{\ref{app:lifecycle}} gives the full lifecycle and recovery details.}

\arxivchange{Keeping the PDHG iterations distributed does not by itself remove the memory limit of one machine. In D-PDLP,\footnote{\arxivchange{The implementation comparison uses D-PDLP commit \href{https://github.com/Lhongpei/D-PDLP/tree/6c1b8174b33665f99ff61abed0539b004519d68f}{\texttt{6c1b817}}.}} one process reads the complete LP and creates a second full copy for scaling before distributing the matrix blocks. At termination, the primal, dual, and reduced-cost vectors are assembled on that process. The input matrix and its scaling copy must therefore fit in one machine's host memory, even when the distributed iteration fits across the available GPUs.}

\arxivchange{\solvername{} instead scales the local matrix blocks, combining row and column statistics across the corresponding process groups. It uses the distributed operator for feasibility polishing as well as optimization, and recovers and writes the solution in slices. No process assembles the complete matrix or solution for these stages. The matrix partition can be constructed from an MPS file during initialization; separately prepared input blocks are optional. Appendix~\ref{app:lifecycle} gives the data layout for each stage.}

\arxivchange{\section{Memory layout and peak memory}
\label{sec:memory}}

\arxivchange{Once the matrix is partitioned, each GPU still needs several vectors for the current iterate, restart anchor, residuals, and temporary calculations. Primal vectors have one entry per local variable, and dual vectors have one entry per local row. A double-precision vector with one billion entries occupies 8~GB; two extra copies cost another 16~GB. Fitting the matrix alone therefore does not ensure that an iteration or a polishing step will fit.}

\arxivchange{We reduce this storage cost through both representation and reuse. Global dimensions and block offsets use 64-bit integers, while indices within a matrix block use 32 bits. Each local index occupies four bytes, and the block offset supplies its position in the complete LP. The local row, column, and nonzero counts must fit the 32-bit representation. \arxivgreen{The complete LP can therefore contain more than $2^{31}-1$ nonzeros, provided each local block fits.} A bound array can also be omitted when all lower bounds are zero, or all upper bounds are positive infinity: the projection uses the constant directly.}

\arxivchange{Temporary vectors need not all coexist. For example, the matrix product reads the reflected primal iterate before that storage is used for an iterate difference and, later, a stationarity residual. Sharing one buffer for these three purposes eliminates two primal-length arrays. \arxivgreen{Figure~\ref{fig:memory-reuse} illustrates these reductions in vector storage.} The distributed Google benchmark runs use this reuse and omit constant bound arrays when the conditions above hold. Appendix~\ref{app:memory} gives the byte counts and explains the additional working vectors needed by feasibility polishing.}

\arxivgreen{\begin{figure}[!ht]
  \centering
  {\color{black}\resizebox{\linewidth}{!}{
\begingroup
\definecolor{memblue}{RGB}{0,114,178}
\definecolor{memorng}{RGB}{210,140,0}
\definecolor{memteal}{RGB}{0,133,119}
\begin{tikzpicture}[
  x=1cm,y=1cm,font=\fontsize{9}{10.5}\selectfont,
  panel/.style={draw=black!23,fill=black!1,rounded corners=2pt,line width=.45pt},
  flow/.style={-{Latex[length=1.8mm,width=1.1mm]},draw=black!65,line width=.55pt},
  outline/.style={draw=black!28,fill=white,line width=.45pt},
  labeltext/.style={align=center,inner sep=1pt}
]
\path[use as bounding box] (0,-.62) rectangle (15.8,7.8);
\draw[panel] (0,0) rectangle (9.85,7.8);
\draw[panel] (10.1,0) rectangle (15.8,7.8);
\node[anchor=west,font=\fontsize{10.5}{12}\selectfont\bfseries] at (.23,7.42)
  {(a) One buffer, three uses};
\node[anchor=west,font=\fontsize{10.5}{12}\selectfont\bfseries] at (10.33,7.42)
  {(b) Omit uniform bounds};

\node[labeltext] at (3.2,6.80) {PDHG products\\and updates};
\node[labeltext] at (5.4,6.80) {Fixed-point\\calculation};
\node[labeltext] at (7.6,6.80) {Stationarity\\check};
\draw[flow] (2.12,6.35) -- (8.68,6.35);

\node[labeltext,anchor=west] at (.22,4.88) {Separate\\arrays};
\foreach \y in {5.70,4.88,4.06} {
  \draw[outline] (2.1,\y-.36) rectangle (8.7,\y+.36);
}
\fill[memblue!16] (2.11,5.35) rectangle (4.3,6.05);
\fill[memorng!20] (4.3,4.53) rectangle (6.5,5.23);
\fill[memteal!17] (6.5,3.71) rectangle (8.69,4.41);
\node[labeltext,text=memblue!85!black] at (3.2,5.70) {reflected\\iterate};
\node[labeltext,text=memorng!70!black] at (5.4,4.88) {iterate\\difference};
\node[labeltext,text=memteal!90!black] at (7.6,4.06) {stationarity\\residual};
\node at (5.4,3.45) {3 allocations: $24n_c$ bytes};

\node[labeltext,anchor=west] at (.22,2.79) {Shared\\array};
\draw[outline] (2.1,2.43) rectangle (8.7,3.15);
\fill[memblue!16] (2.11,2.44) rectangle (4.3,3.14);
\fill[memorng!20] (4.3,2.44) rectangle (6.5,3.14);
\fill[memteal!17] (6.5,2.44) rectangle (8.69,3.14);
\draw[white,line width=1pt] (4.3,2.44) -- (4.3,3.14);
\draw[white,line width=1pt] (6.5,2.44) -- (6.5,3.14);
\node[labeltext,text=memblue!85!black] at (3.2,2.79) {reflected\\iterate};
\node[labeltext,text=memorng!70!black] at (5.4,2.79) {iterate\\difference};
\node[labeltext,text=memteal!90!black] at (7.6,2.79) {stationarity\\residual};
\node at (5.4,2.08) {1 allocation: $8n_c$ bytes};
\draw[flow] (9.10,4.88) -- (9.10,2.79);

\node[labeltext,text width=5.0cm] at (12.95,6.66)
  {If every local variable has bounds $[0,+\infty)$};
\draw[outline] (10.48,5.55) rectangle (15.42,6.10);
\draw[outline] (10.48,4.71) rectangle (15.42,5.26);
\node[anchor=west] at (10.66,5.82) {$\ell$: $0,\ 0,\ 0,\ \ldots,\ 0$};
\node[anchor=east] at (15.25,5.82) {$8n_c$};
\node[anchor=west] at (10.66,4.98) {$u$: $\infty,\ \infty,\ \ldots,\ \infty$};
\node[anchor=east] at (15.25,4.98) {$8n_c$};
\node at (12.95,4.20) {Two arrays: $16n_c$ bytes};
\draw[flow] (12.95,3.90) -- (12.95,3.35);
\draw[draw=memteal,fill=memteal!7,rounded corners=2pt]
  (10.48,2.52) rectangle (15.42,3.22);
\node[labeltext] at (12.95,2.87)
  {Projection: $v_j\mapsto\max\{v_j,0\}$};
\node[labeltext] at (12.95,2.06) {No per-variable bound arrays};

\draw[draw=none,fill=memteal!8,rounded corners=2pt]
  (.24,.26) rectangle (9.61,1.53);
\node[font=\bfseries,text=memteal!85!black] at (4.925,1.20)
  {$16n_c$ bytes saved};
\node[labeltext] at (4.925,.78)
  {For $n_c=10^9$: 24 GB $\longrightarrow$ 8 GB};
\node[font=\scriptsize,text=black!65] at (4.925,.43)
  {FP64 storage for these three vectors only};
\draw[draw=none,fill=memteal!8,rounded corners=2pt]
  (10.34,.26) rectangle (15.56,1.53);
\node[font=\bfseries,text=memteal!85!black] at (12.95,1.20)
  {$16n_c$ bytes saved};
\node[labeltext] at (12.95,.76)
  {For $n_c=10^9$: 16 GB saved};
\node[font=\scriptsize,text=black!65] at (12.95,.43)
  {when both bounds qualify};
\node[font=\fontsize{8}{9}\selectfont,text=black!65] at (7.9,-.34)
  {$n_c$: number of local variables. Each FP64 entry uses 8 bytes; GB $=10^9$ bytes.};
\end{tikzpicture}
\endgroup}}
  \caption{\arxivgreen{Reducing vector storage on one GPU. (a) Values used at different times share one buffer, saving two FP64 vectors. Colored regions show when each value is used. (b) Uniform zero lower bounds and infinite upper bounds need no arrays. Here $n_c$ is the number of local variables; the byte counts cover only the illustrated arrays.}}
  \label{fig:memory-reuse}
\end{figure}}
\arxivgreen{\FloatBarrier}

\section{Computational Results}

We organize the experiments around three questions: \textbf{RQ1}: Can \solvername{} solve the Google PDLP benchmark? \textbf{RQ2}: How far can it scale beyond that benchmark? \textbf{RQ3}: Does communication limit multi-GPU scaling, and when does support-aware communication help? All experiments ran on Georgia Tech's PACE Phoenix Cluster, with one MPI rank per GPU. A compute node is one physical server. The validated runs use 1--76 GPUs and up to 29 nodes; Appendix~\ref{app:hardware} gives the hardware, interconnect, software, and memory details.
\paragraph{Validation and timing.}
\arxivremovedparagraph{For RQ1, we use the Google study's published scaled-LP criterion \cite{applegate2026pdlp}. For every other reported solve, a separate checker evaluates the exported primal--dual solution against the original unscaled LP. Acceptance requires all nine quantities in Equation~\mbox{\eqref{eq:uniform-check}}, covering normalized primal feasibility, stationarity, dual-sign admissibility, and relative primal--dual gap, to be finite and at most $10^{-6}$. Solver time covers optimization and feasibility polishing when enabled; end-to-end time additionally includes sharded input, initialization, recovery, and output. Offline shard construction and the final independent check are excluded. Appendix~\mbox{\ref{app:protocol}} gives the full protocol.}

\arxivchange{For the solved RQ1 instances other than Heat Source Easy, we use the Google study's published scaled-LP criterion \cite{applegate2026pdlp}. Heat Source Easy and the solves in RQ2 and RQ3 are checked against the original unscaled LP. All nine quantities in Equation~\eqref{eq:uniform-check} must be finite and no greater than $10^{-6}$, covering primal feasibility, stationarity, dual-sign admissibility, and relative primal--dual gap. This independent check is applied to the exported vectors after the solver has terminated under its internal stopping rule.}

\arxivchange{Solver time \arxivgreen{in RQ1 and RQ2} measures elapsed optimization through termination, including convergence evaluations, restarts, and all feasibility-polishing attempts. End-to-end time also includes loading the prepared matrix blocks, initialization, solution recovery, and output. All reported runs use prepared input; constructing those blocks and performing the final independent check are excluded from both clocks. Appendix~\ref{app:protocol} gives the stopping rules and timing details. \arxivgreen{The paired communication comparisons in RQ3 use the timing definitions in Appendix~\ref{app:hardware}.}}

\subsection{RQ1: Google PDLP benchmark}

\begin{table}[H]
\centering
\scriptsize
\setlength{\tabcolsep}{2.05pt}
\renewcommand{\arraystretch}{1.13}
\begin{tabular}{@{}l r r r c r r r r@{}}
\toprule
& \multicolumn{3}{c}{Problem size}
& \multicolumn{3}{c}{\cellcolor{shardfill}\textcolor{sharddark}{\textbf{\solvername{} -- SECONDS}}}
& \multicolumn{2}{c}{\cellcolor{cpufill}\textcolor{cpudark}{\textbf{Published CPU -- HOURS}}} \\
\cmidrule(lr){2-4}\cmidrule(lr){5-7}\cmidrule(l){8-9}
Instance
& $m$ & $n$ & $\nnz(A)$
& \cellcolor{shardfill}GPUs/nodes
& \cellcolor{shardfill}solver
& \cellcolor{shardfill}end-to-end
& \cellcolor{cpufill}Google PDLP
& \cellcolor{cpufill}Gurobi \\
\midrule
TSP-Gaia-100M
& 162.935M & 1.185B   & 6.338B
& \cellcolor{shardfill}8/1
& \cellcolor{shardfill}\textbf{593\,s}
& \cellcolor{shardfill}933\,s
& \cellcolor{cpufill}21.06\,h
& \cellcolor{cpufill}-- \\
Design Match
& 22.000M  & 40.000M  & 2.760B
& \cellcolor{shardfill}2/1
& \cellcolor{shardfill}\textbf{232\,s}
& \cellcolor{shardfill}373\,s
& \cellcolor{cpufill}9.34\,h
& \cellcolor{cpufill}32.4\,h \\
QAP-THO-150
& 6.705M   & 249.784M & 1.006B
& \cellcolor{shardfill}4/2
& \cellcolor{shardfill}\textbf{133\,s}
& \cellcolor{shardfill}160\,s
& \cellcolor{cpufill}--
& \cellcolor{cpufill}-- \\
World Shipping
& 15.304M  & 228.868M & 688.659M
& \cellcolor{shardfill}2/1
& \cellcolor{shardfill}\textbf{1,921\,s}
& \cellcolor{shardfill}2,000\,s
& \cellcolor{cpufill}53.81\,h
& \cellcolor{cpufill}-- \\
Mediterranean
& 7.491M   & 208.479M & 628.927M
& \cellcolor{shardfill}2/1
& \cellcolor{shardfill}\textbf{5,839\,s}
& \cellcolor{shardfill}5,890\,s
& \cellcolor{cpufill}32.62\,h
& \cellcolor{cpufill}-- \\
Production Inventory
& 4.651M   & 18.271M  & 500.050M
& \cellcolor{shardfill}1/1
& \cellcolor{shardfill}\textbf{2,305\,s}
& \cellcolor{shardfill}2,358\,s
& \cellcolor{cpufill}--
& \cellcolor{cpufill}5.8\,h \\
TSP-Gaia-10M
& 17.017M  & 60.602M  & 475.702M
& \cellcolor{shardfill}1/1
& \cellcolor{shardfill}\textbf{243\,s}
& \cellcolor{shardfill}335\,s
& \cellcolor{cpufill}2.99\,h
& \cellcolor{cpufill}28.1\,h \\
Supply Chain
& 2.210M   & 201.000M & 403.000M
& \cellcolor{shardfill}1/1
& \cellcolor{shardfill}\textbf{162\,s}
& \cellcolor{shardfill}224\,s
& \cellcolor{cpufill}18.90\,h
& \cellcolor{cpufill}2.5\,h \\
QAP-WIL-100
& 1.980M   & 49.015M  & 198.020M
& \cellcolor{shardfill}2/1
& \cellcolor{shardfill}\textbf{15.9\,s}
& \cellcolor{shardfill}44.6\,s
& \cellcolor{cpufill}0.28\,h
& \cellcolor{cpufill}-- \\
Heat Source Easy\arxivchange{$^{*}$}
& 15.625M  & 31.628M  & 125.000M
& \cellcolor{shardfill}\arxivremoved{--}\,\arxivchange{8/8}
& \cellcolor{shardfill}\arxivchange{\textbf{28,858\,s}}
& \cellcolor{shardfill}\arxivchange{28,878\,s}
& \cellcolor{cpufill}59.97\,h
& \cellcolor{cpufill}-- \\
Heat Source Hard
& 15.625M  & 31.628M  & 125.000M
& \multicolumn{3}{c}{\cellcolor{shardfill}--}
& \cellcolor{cpufill}--
& \cellcolor{cpufill}-- \\
\bottomrule
\end{tabular}
\caption{Google PDLP benchmark results \citep{applegate2026pdlp}. \textcolor{sharddark}{\textbf{\solvername{} runtimes are reported in seconds}}; \textcolor{cpudark}{\textbf{published CPU runtimes are reported in hours}}. A dash marks a target not reached. Instances are ordered by decreasing $\nnz(A)$. \arxivchange{$^{*}$Heat Source Easy uses Quadro RTX~6000 GPUs and the original-space $10^{-6}$ checks in Equation~\eqref{eq:uniform-check}; the other GPU rows use H200 GPUs and the published Google criterion.}}
\label{tab:google}
\end{table}
 

\solvername{} reaches the published PDLP criterion on nine of the eleven benchmark instances, versus eight in the published CPU PDLP study (Table~\ref{tab:google}). \arxivchange{Heat Source Easy adds a tenth solved instance under the original-space $10^{-6}$ checks.} On TSP-Gaia-100M, eight H200 GPUs require 592.6 seconds of solver time and 933.4 seconds end to end; the published 32-core CPU experiment reports 21.06 hours.

\subsection{RQ2: Scaling beyond the Google benchmark}

\arxivgreenreplace{Table~\mbox{\ref{tab:scale}} pushes the same distributed solve beyond the Google benchmark.}{Table~\ref{tab:scale} reports validated solves beyond the Google benchmark.} It includes KDD12 \citep{chang2011libsvm,libsvmdata2026}, QAPLIB \citep{adams1994improved,burkard1997qaplib} and MS1 matching \citep{koohi2024msbiographs} cases, plus a deterministic multicommodity-flow family used for the largest scale tests. MCF5B is validated in all three runs on 40 GPUs across ten compute nodes. Larger instances extend the results to 25.000 billion nonzeros on 56 GPUs and 40.807 billion nonzeros on 76 GPUs across 29 nodes. \arxivgreenreplace{The latter has 13.604 billion variables}{MCF13.60B has 13.604 billion variables} and solves in 1,905 seconds. The KDD12, QAPLIB and MS1 cases show that the distributed path is not limited to the synthetic flow family. Within the comparison scope in Appendix~\ref{app:scope}, MCF13.60B is, to our knowledge, the largest reported GPU LP solve with a separately validated complete primal--dual solution.

\begin{table}[!t]
\centering
\scriptsize
\setlength{\tabcolsep}{2.5pt}
\renewcommand{\arraystretch}{1.13}
\begin{tabular*}{\linewidth}{@{\extracolsep{\fill}}llrrrcrr@{}}
\toprule
& & \multicolumn{3}{c}{Problem size}
& \multicolumn{3}{c}{\cellcolor{shardfill}\textcolor{sharddark}{\textbf{\solvername{} -- SECONDS}}} \\
\cmidrule(lr){3-5}\cmidrule(l){6-8}
Instance
& source
& $m$
& $n$
& $\nnz(A)$
& \cellcolor{shardfill}GPUs/nodes
& \cellcolor{shardfill}solver
& \cellcolor{shardfill}end-to-end \\
\midrule
MCF13.60B & synthetic & 2.040M & 13.604B & 40.807B
& \cellcolor{shardfill}76 mixed / 29
& \cellcolor{shardfill}\textbf{1,905\,s}
& \cellcolor{shardfill}2,057\,s \\
MCF8.35B & synthetic & 24.951M & 8.353B & 25.000B
& \cellcolor{shardfill}56 mixed / 25
& \cellcolor{shardfill}\textbf{2,933\,s}
& \cellcolor{shardfill}3,045\,s \\
MCF8B
& synthetic
& 101.059M
& 8.000B
& 16.042B
& \cellcolor{shardfill}12 H200 / 3
& \cellcolor{shardfill}\textbf{14,905\,s}
& \cellcolor{shardfill}15,510\,s \\
MCF5B
& synthetic
& 63.162M
& 5.000B
& 10.027B
& \cellcolor{shardfill}40 RTX / 10
& \cellcolor{shardfill}\textbf{560\,s}
& \cellcolor{shardfill}634\,s \\
\texttt{tai256c} AJ
& QAPLIB
& 33.424M
& 2.131B
& 8.557B
& \cellcolor{shardfill}8 BW / 1
& \cellcolor{shardfill}\textbf{717\,s}
& \cellcolor{shardfill}920\,s \\
KDD12 L1-SVM
& public data
& 149.639M
& 259.012M
& 3.442B
& \cellcolor{shardfill}8 H200 / 1
& \cellcolor{shardfill}\textbf{37,294\,s}
& \cellcolor{shardfill}37,587\,s \\
MS1 matching
& public data
& 43.144M
& 1.309B
& 2.617B
& \cellcolor{shardfill}32 V100 / 16
& \cellcolor{shardfill}\textbf{3,721\,s}
& \cellcolor{shardfill}3,760\,s \\
\bottomrule
\end{tabular*}
\caption{Validated scale experiments beyond the Google benchmark. All rows are separately validated solves and are ordered by decreasing $\nnz(A)$. BW and RTX denote RTX PRO 6000 Blackwell and Quadro RTX~6000, respectively. Mixed GPU allocations are detailed in Appendix~\ref{app:hardware}.}
\label{tab:scale}
\end{table}

\FloatBarrier

\subsection{RQ3: Does communication limit multi-GPU scaling?}

\begin{figure}[h!]
\centering
\includegraphics[width=\linewidth]{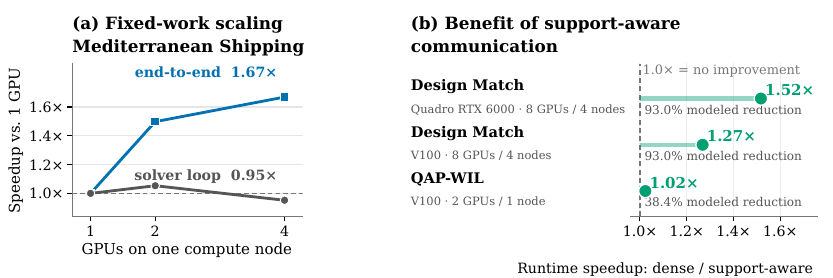}
\caption{Communication effects. (a) Speedup for a fixed 4,000 Mediterranean Shipping iterations relative to one GPU. (b) Geometric-mean speedup of support-aware over dense communication; $1\times$ means no change.}
\label{fig:med-participant}
\end{figure}

Adding GPUs reduces local matrix work but adds communication. On Mediterranean Shipping, 4,000 iterations take 107.2 seconds on one RTX PRO 6000 Blackwell GPU and 112.7 seconds on four GPUs: the solver loop is slightly slower, even though end-to-end time improves from 281.5 to 168.1 seconds (Figure~\ref{fig:med-participant}a). Thus, smaller per-GPU matrix blocks do not automatically make the solver loop faster.

With a column partition, each GPU stores only some columns of $A$. A dense reduction still involves every GPU when forming $Ax$, even when a GPU has no nonzeros in a row. Support-aware communication skips those GPUs for that row. Because their contribution is zero, the exact-arithmetic PDHG update is unchanged; Appendix~\ref{app:participant} gives the details.

Figure~\ref{fig:med-participant}b shows that the benefit depends on \arxivgreenreplace{how many GPUs actually contain each row}{how many GPUs store nonzeros in each row}. On Design Match, modeled communication falls by 92.97\%, and support-aware communication is $1.269\times$ faster on V100s and $1.518\times$ faster on Quadro RTX~6000s. On QAP-WIL, the modeled reduction is only 38.35\% and the speedup is $1.024\times$. Every paired run passes the separate original-space $10^{-6}$ check. Pairing and timing details are in Appendix~\ref{app:hardware}.

\section{Limitations and Future Work}
\label{sec:limitations}

\paragraph{Limitations.}
\arxivreplaced{\mbox{\solvername{}} assumes that the LP has already been partitioned into shards; converting a monolithic model to this format is currently an offline preprocessing step. In addition, support-aware communication is implemented only for column partitions with fixed matrix support. These choices are sufficient for the experiments in this paper, but they leave room for a more general distributed solver pipeline.}{Input construction and independent validation can still be limited by host memory on one machine. Building blocks from an MPS file requires temporary coefficient storage, while checking against that file keeps complete solution and work vectors in one CPU process. These stages can therefore exceed host memory even when the matrix and iterates fit across the GPUs. Support-aware communication is implemented only for column partitions whose matrix support remains fixed; extending it to general two-dimensional grids requires a different communication plan.}

\paragraph{Future work.}
Future work will focus on scaling to still larger LPs and reducing both communication and memory use. This includes better partitioning, broader distributed presolve and postsolve, and extending support-aware communication to more general process grids. \arxivgreenreplace{We are also interested in applying the same ideas to other large-scale optimization problems, including quadratic and semidefinite programs, and in using \mbox{\solvername{}} as a subroutine inside algorithms for very large discrete optimization problems, such as decomposition methods or large-scale heuristics.}{Further directions include quadratic and semidefinite optimization and the use of \solvername{} within decomposition methods and heuristics for discrete optimization.}

\section{Conclusion}
 
Previous work distributes the main PDHG computations across GPUs. \solvername{} goes further by keeping the problem and solution distributed throughout the solve. Persistent ownership enables solutions \arxivreplaced{meeting the published criterion on nine of eleven Google benchmarks}{for ten of eleven Google benchmarks} and separately checked LPs with 13.604 billion variables and 40.807 billion nonzeros, with validated runs spanning up to 29 compute nodes. Smaller per-GPU blocks alone do not ensure speed: communication must also follow matrix sparsity. On Design Match, support-aware communication cuts modeled data movement by 92.97\% and yields $1.27\times$--$1.52\times$ solver-time speedups. \arxivgreenreplace{Together, these results isolate two distinct scale bottlenecks: memory ownership and communication.}{The results show why memory distribution and communication must be considered together.}

\clearpage
\bibliography{references}

\appendix

\section{Computational environment and hardware}
\label{app:hardware}

All experiments were executed on Georgia Tech's PACE Phoenix cluster under Slurm, with one MPI process per GPU. Table~\ref{tab:hardware} summarizes the node classes used, using retained GPU inventories and the corresponding Slurm node configurations. \arxivreplaced{All solved rows in Table~\mbox{\ref{tab:google}} use H200 GPUs;}{In Table~\ref{tab:google}, Heat Source Easy uses Quadro RTX~6000 GPUs and the other solved rows use H200 GPUs;} Table~\ref{tab:scale} identifies the GPU families for the scale experiments. The GPU/node counts in these tables describe GPUs used by the solver, which may be fewer than the GPUs installed in an allocated node.

\begin{table}[!ht]
\centering
\scriptsize
\setlength{\tabcolsep}{4pt}
\renewcommand{\arraystretch}{1.15}
\caption{PACE Phoenix node classes used in the reported experiments. GPU memory is the capacity reported by \texttt{nvidia-smi}; host RAM is Slurm's configured node memory, rounded to the nearest GiB. CPU counts are configured cores per node. These are node capacities, not per-job reservations.}
\label{tab:hardware}
\begin{tabular}{@{}>{\raggedright\arraybackslash}p{3.3cm}rr>{\raggedright\arraybackslash}p{3.7cm}r@{}}
\toprule
GPU & \shortstack{memory/GPU\\(GiB)} & \shortstack{installed\\GPUs/node}
& Host CPU (cores/node) & \shortstack{host RAM\\(GiB)} \\
\midrule
NVIDIA H200 & 140.4 & 8
& Intel Xeon Platinum 8562Y+ (64) & 2,015 \\
RTX PRO 6000 Blackwell Server Edition & 95.6 & 8
& Intel Xeon, Granite Rapids (64) & 2,015 \\
Quadro RTX 6000 & 24.0 & 4
& Intel Xeon Gold 6226 (24) & 376 / 754 \\
Tesla V100 PCIe & 16.0 / 32.0 & 2
& Intel Xeon Gold 6226 (24) & 376 / 754 \\
\bottomrule
\end{tabular}
\end{table}

\paragraph{Software and communication.}
The distributed C/CUDA runs use CUDA~12.9.1 and CUDA-aware OpenMPI~4.1.8. GPU communication uses NVIDIA Collective Communications Library (NCCL)~2.26.5, with builds targeting each GPU architecture. MPI manages the process grid and control reductions; NCCL exchanges GPU buffers. Retained \texttt{nvidia-smi topo -m} output shows NVLink connectivity between the H200 GPUs (\texttt{NV18}), and PCIe paths on the Blackwell, Quadro RTX~6000, and V100 nodes used here.

Inter-node transport is run-specific. MCF8B, MCF5B, MS1, and the Design Match communication comparisons disable NCCL's InfiniBand transport (\texttt{NCCL\_IB\_DISABLE=1}) and use its socket path. The two-node QAP-THO-150 benchmark run enables the InfiniBand transport. \arxivgreenreplace{These settings do not establish a common physical-link bandwidth or a cluster-wide RDMA performance claim.}{Because transport settings differ across experiments, the timings do not isolate the effect of RDMA.}

\paragraph{Allocation and timing conditions.}

MCF8.35B combines 42 Quadro RTX~6000 and 14 V100-32GB GPUs across 25 nodes. MCF13.60B combines 60 Quadro RTX~6000, 14 V100-32GB, and two L40S GPUs across 29 nodes. Both use column partitions, one MPI rank per GPU, and socket transport.

Jobs reserve GPUs, CPU cores, and host memory through Slurm; whole-node exclusivity is not assumed. For example, MCF8B uses four of the eight H200s on each of three nodes, and the Quadro RTX~6000 Design Match comparisons use two of four GPUs per node. The Mediterranean 1/2/4-GPU fixed-work measurements use exclusive Blackwell nodes, with all three GPU counts tested within each allocation. Dense/support-aware pairs run sequentially on the same allocation and GPU set, with order reversed across pairs. Figure~\ref{fig:med-participant}b reports geometric means over two order-balanced pairs per Design Match platform and eight QAP-WIL pairs across four allocations. For Design Match, the timing metric is solver time through the first validated checkpoint; for QAP-WIL, it is elapsed time through validation. MCF4.15B uses two shared nodes; its timing is descriptive rather than a scaling measurement.

\begin{table}[!ht]
\centering
\scriptsize
\setlength{\tabcolsep}{4pt}
\caption{Physical footprint of selected checked executions. Sizes are GiB ($2^{30}$ bytes); input and primal--dual output are totals across all shards. Peak GPU memory is the largest sampled device-memory use on any participating GPU. BW denotes RTX PRO 6000 Blackwell Server Edition.}
\label{tab:storage-memory}
\begin{tabular}{@{}lccrrrr@{}}
\toprule
Instance & GPU & process grid & GPUs/nodes & input & output &
\shortstack{peak GPU\\memory/GPU} \\
\midrule
TSP-Gaia-100M & H200 & $1\times8$ & 8/1 & 121.59 & 18.87 & 80.95 \\
\texttt{tai256c} AJ & BW & $1\times8$ & 8/1 & 148.23 & 32.00 & 69.39 \\
MCF8B & H200 & $1\times12$ & 12/3 & 380.69 & 119.96 & 79.78 \\
\bottomrule
\end{tabular}
\end{table}

Table~\ref{tab:storage-memory} reports physical storage and device-memory use. Peak memory is the maximum of periodic \texttt{nvidia-smi memory.used} samples, matched to the participating GPU UUIDs across the allocation. \arxivgreenreplace{It is a sampled device-level quantity, not an exact allocator high-water mark.}{Short-lived allocation peaks may fall between samples.}


\section{Instance origins}
\label{app:instances}

Here $m$, $n$, and $\nnz(A)$ denote constraints, variables, and explicitly stored matrix coefficients. Table~\ref{tab:google-provenance} summarizes the released dimensions and construction of the Google PDLP benchmark \citep{applegate2026pdlp}.\footnote{The scaled instance files are available from Oliver Hinder's \href{https://www.oliverhinder.com/large-scale-lp-problems}{benchmark download page}; generators and detailed construction notes for the subset with released construction code are available in the companion \href{https://github.com/ohinder/large-scale-LP-test-problems}{GitHub repository}.}  We use \emph{synthetic} for generated, application-inspired models and identify LPs constructed from a named public data source separately.

\begin{table}[!ht]
\centering
\caption{Meaning and construction of the eleven Google benchmark instances.}
\label{tab:google-provenance}
\scriptsize
\setlength{\tabcolsep}{2.0pt}
\begin{tabular}{@{}lrrr>{\raggedright\arraybackslash}p{6.0cm}@{}}
\toprule
Instance & $m$ & $n$ & $\nnz(A)$ & origin and model \\
\midrule
Design Match & 22,000,135 & 40,000,000 & 2,760,000,000
& Synthetic covariate-balancing statistical-matching LP from the application
class studied by Zubizarreta
\citep{zubizarreta2012matching}. \\
Heat Source Easy & 15,625,000 & 31,628,008 & 125,000,000
& Synthetic inverse heat-source problem with temperature observations and
linear PDE constraints. \\
Heat Source Hard & 15,625,000 & 31,628,008 & 125,000,000
& The same synthetic inverse problem with fewer measurements and more true and
candidate source locations. \\
Mediterranean Shipping & 7,490,593 & 208,479,461 & 628,927,462
& LP relaxation of a liner-shipping mixed-integer program (MIP) built from the LINERLIB Mediterranean
benchmark \citep{brouer2014linerlib}. \\
Production Inventory & 4,650,850 & 18,270,600 & 500,049,700
& Synthetic robust production--inventory model based on Ben-Tal et al.
\citep{bental2004adjustable}, with randomized data. \\
QAP--THO--150 & 6,705,300 & 249,783,750 & 1,005,795,000
& Adams--Johnson LP relaxation \citep{adams1994improved} of the QAPLIB
\texttt{tho150} benchmark \citep{burkard1997qaplib}. \\
QAP--WIL--100 & 1,980,200 & 49,015,000 & 198,020,000
& Adams--Johnson LP relaxation \citep{adams1994improved} of the QAPLIB
\texttt{wil100} benchmark \citep{burkard1997qaplib}. \\
Supply Chain & 2,210,100 & 201,000,100 & 403,000,100
& Synthetic multicommodity-flow model for large-retailer supply-chain planning. \\
TSP--Gaia--100M & 162,934,799 & 1,184,557,727 & 6,337,834,450
& Public-data-derived TSP lower-bound LP on the 100 million nearest stars in
Gaia DR2 \citep{brown2018gaia,cook2024gaia100m}; degree constraints plus
62,934,799 cuts collected with Concorde \citep{applegate2020concorde}. \\
TSP--Gaia--10M & 17,016,681 & 60,601,996 & 475,701,996
& Public-data-derived counterpart on the ten million nearest stars in Gaia DR2
\citep{brown2018gaia,cook2024gaia10m}; degree constraints plus 7,016,681 cuts
collected with Concorde \citep{applegate2020concorde}. \\
World Shipping & 15,304,282 & 228,867,510 & 688,658,522
& LP relaxation of a liner-shipping MIP built from the LINERLIB World
benchmark \citep{brouer2014linerlib}. \\
\bottomrule
\end{tabular}
\end{table}

\begin{table}[!ht]
\centering
\caption{Origins and validation status of the beyond-benchmark
experiments.}
\label{tab:beyond-provenance}
\scriptsize
\setlength{\tabcolsep}{2.0pt}
\begin{tabular}{@{}lrrr>{\raggedright\arraybackslash}p{6.0cm}@{}}
\toprule
Instance & $m$ & $n$ & $\nnz(A)$ & origin, model, and status \\
\midrule
KDD12 L1-SVM & 149,639,105 & 259,012,009 & 3,441,699,415
& Public-data-derived LP from the LIBSVM \texttt{kdd12.xz} click-log data,
using a no-intercept, $C=1$, split-variable $\ell_1$-regularized hinge-loss model
\citep{chang2011libsvm,libsvmdata2026}; one validated run. \\
MS1 matching & 43,144,218 & 1,308,742,322 & 2,617,484,644
& Maximum-weight fractional matching LP from the MS-BioGraphs MS1 protein
sequence-similarity graph \citep{koohi2024msbiographs}; one variable per
undirected non-loop edge, unit vertex capacities, and $0\le x_e\le1$.
Self-loops and reverse duplicates are removed; one original-LP checked run. \\
\texttt{tai256c} AJ & 33,423,872 & 2,130,804,736 & 8,556,511,232
& Benchmark-derived, symmetry-reduced Adams--Johnson/substitution LP generated
from QAPLIB \texttt{tai256c}
\citep{adams1994improved,burkard1997qaplib}; one separately checked execution. \\
MCF4.15B & 52,423,770 & 4,149,986,400 & 8,321,930,400
& Deterministic synthetic multicommodity-flow model; one separately checked
execution on eight RTX PRO 6000 Blackwell GPUs across two shared compute nodes,
used only to demonstrate feasibility, not performance or scaling. \\
MCF5B & 63,161,790 & 5,000,032,800 & 10,026,520,800
& Deterministic synthetic multicommodity-flow model with commodity-local
factory--warehouse--store blocks; three separately checked executions. \\
MCF8B & 101,059,055 & 8,000,067,600 & 16,042,463,600
& Larger member of the same deterministic synthetic family; two separately
checked executions. \\
MCF8.35B & 24,950,515 & 8,352,818,000 & 25,000,333,000
& Same deterministic MCF family; one checked execution on
56 GPUs across 25 nodes. \\
MCF13.60B & 2,040,170 & 13,604,080,000 & 40,807,480,000
& Same family with larger commodity blocks; one checked execution
on 76 GPUs across 29 nodes. \\
MCF8.5B capacity & 107,374,470 & 8,500,010,400 & 17,044,994,400
& Capacity-only run from the same family on 70 GPUs; complete solution vectors
were not exported, so the run is not counted as a solve. \\
\bottomrule
\end{tabular}
\end{table}

MCF8.35B uses 8,303 commodities, 1,000 factories per commodity, 1,000 warehouses, and five stores. MCF13.60B uses 34 commodities, 20,000 factories per commodity, 20,000 warehouses, and five stores. The latter has more nonzeros but fewer rows, reducing row-vector communication in the column-partitioned solve.

\FloatBarrier

\section{Distributed solver stages and solution validation}
\label{app:protocol}

\subsection{Data layout across solver stages}
\label{app:lifecycle}

\arxivreplaced{The solve begins from prepared matrix shards. Conversion from a monolithic input (e.g., an MPS file) to these shards is performed offline and is not included in reported times.}{Initialization can construct the matrix blocks from an MPS file, or each process can read a block prepared in advance. Both routes produce the same distributed layout for optimization.} In an $R\times C$ grid, each process holds $A_{rc}$, $x_c$, and $y_r$, with vector replicas only within the corresponding grid column or row. Table~\ref{tab:lifecycle} specifies how this ownership is retained.

\begin{table}[!ht]
\centering
\caption{Implemented distributed lifecycle.  ``Global reduction'' means a
scalar or partition-sized collective, never a complete matrix or primal--dual
gather. \arxivchange{Independent validation is performed separately, as described below.}}
\label{tab:lifecycle}
\scriptsize
\setlength{\tabcolsep}{3.0pt}
\begin{tabular}{@{}>{\raggedright\arraybackslash}p{2.25cm}p{9.85cm}@{}}
\toprule
Phase & Ownership and communication \\
\midrule
\arxivreplaced{Prepared-shard input}{Input} & \arxivreplaced{Each process reads its matrix block and matching row/column data directly.  Offline construction and conversion are excluded; no claim is made for arbitrary monolithic-input conversion.}{Initialization constructs each process's matrix block from the MPS input, or loads that block and its row/column data from prepared files.} \\
Scaling & Local nonzeros form row and column statistics; grid reductions form the diagonal scales, which remain partitioned with the vectors. \\
Optimization & Local products plus the two reductions in Equation~\eqref{eq:distributed-products}; Halpern, projection, and vector updates act on local blocks. \\
Restart/termination & Processes reduce norms and scalar decision statistics; latest, average, anchor, and reflected states retain the same block layout. \\
Singleton-zero presolve & One exact elimination pass identifies zero-equality singleton rows, fixes their columns at zero, removes those entries locally, and stores a distributed pivot map.  This is not arbitrary general presolve. \\
Feasibility polishing & Primal and dual phases reuse the distributed operator; candidate selection and restart decisions use collective statistics while phase states remain local. \\
\arxivreplaced{Recovery, output, and check}{Recovery and output} & The pivot map reconstructs eliminated solution coordinates on their replicas; each process writes owned primal, dual, and reduced-cost slices.  \arxivremoved{The separate checker consumes those shards without a root-vector gather.} \\
\bottomrule
\end{tabular}
\end{table}

\arxivchange{The reported experiments do not use general presolve. Both variants in the Design Match communication comparison apply the same singleton-zero elimination described in Table~\ref{tab:lifecycle}. This removes 1.123 billion nonzeros (40.69\%), leaving 1.637 billion. The Google benchmark runs do not use this reduction.}

\noindent \arxivgreenreplace{For a singleton pivot $(j,s,a_{sj})$, recovery forms}{For an eliminated variable $j$, let $s$ be its singleton zero-equality row and $a_{sj}$ its nonzero coefficient. Recovery forms} the original-unit reduced cost $r=c-A^\top y$ and applies $x_j=0$, $y_s\leftarrow y_s+r_j/a_{sj}$, and $r_j=0$.  This recovery covers primal--dual solutions; infeasibility and unboundedness rays are not implemented.  Prepared blocks use contiguous row and column intervals in a fixed coordinate order; the partition does not change during a solve.

\arxivchange{Input construction and independent validation have different memory requirements from the optimization itself. Reading an MPS file requires tables that map row and column names to coordinates, together with temporary coefficient storage while the matrix blocks are built. The checker for this input format reads matrix coefficients sequentially, but keeps complete primal, dual, reduced-cost, and work vectors in one CPU process. Its host-memory requirement therefore grows as $O(m+n)$ even though it does not retain the complete matrix. The checker that reads prepared blocks instead distributes these vectors. The MPS reader also uses signed 32-bit global row indices before partitioning, which limits the total number of rows that can be imported from MPS. Neither checker is part of the timed optimization.}

\FloatBarrier

\subsection{Validation criteria and timing}

For Table~\ref{tab:google}, general presolve and the singleton-zero reduction are disabled\arxivremoved{, and the scaling from the Google PDLP study is used}. \arxivchange{The solved benchmark instances other than Heat Source Easy use the scaling and acceptance criterion from the Google PDLP study.} \arxivremoved{Both arms of the later Design Match communication comparisons use the same singleton-zero reduction.}  The published rule requires absolute $\ell_\infty$ primal and stationarity residuals at most $10^{-8}$ on the scaled LP, primal/dual bound and sign membership, and a relative primal--dual gap at most $10^{-2}$.  Our finite-precision membership tests use a $10^{-8}$ tolerance \citep{applegate2026pdlp}.

For \arxivreplaced{every other run counted as solved}{Heat Source Easy and the solves in RQ2 and RQ3}, we use a uniform acceptance tolerance of $10^{-6}$.  The separate checker reads the postsolved, sharded solution in the coordinates of the original unscaled LP. Its inputs are the primal vector $x$, the row-dual vector $y$, and the exported reduced-cost vector $r$. It tests four requirements: primal feasibility, stationarity, admissible signs for $y$ and $r$, and agreement between the primal and dual objectives.

\paragraph{Interval quantities.}
For an interval $I=[\ell,u]$, define
\begin{gather*}
 v_I(t):=\operatorname{dist}(t,I),\qquad
 s(I):=\max\bigl(\{0\}\cup
 \{|b|:b\in\{\ell,u\},\ b\text{ finite}\}\bigr),\\
 D(I):=\{q\in\mathbb{R}:
 q>0\Rightarrow\ell>-\infty,\quad
 q<0\Rightarrow u<+\infty\},\qquad
 \psi_I(q):=\inf_{t\in I}qt .
\end{gather*}
Here, $v_I(t)$ is the violation of the interval, and $s(I)$ is a scale
computed from its finite endpoints.  The set $D(I)$ gives the multiplier signs
allowed by the interval: a positive multiplier requires a finite lower bound,
and a negative multiplier requires a finite upper bound.  The function
$\psi_I(q)$ is the contribution of that interval to the dual objective.

\paragraph{Stationarity and objective values.}
We form sign-admissible copies of the reported dual quantities,
\[
 \bar y_i=\proj_{D(S_i)}(y_i),
 \qquad
 \bar r_j=\proj_{D(X_j)}(r_j).
\]
These projections are used only when evaluating stationarity and the dual objective.  \arxivgreenreplace{They do not hide invalid reported signs: the distances between $(y,r)$ and $(\bar y,\bar r)$ are checked separately below.}{The distances between $(y,r)$ and $(\bar y,\bar r)$ are checked separately below to measure violations of the required signs.}

Using the projected quantities, define the stationarity error
\[
 e=c-A^\top\bar y-\bar r,
\]
the primal objective
\[
 p=c_0+c^\top x,
\]
and the dual objective
\[
 d=c_0+\sum_i\psi_{S_i}(\bar y_i)
       +\sum_j\psi_{X_j}(\bar r_j),
\]
where $c_0=0$ if the model has no objective constant.

For each variable and row, define the primal violations and their data scales:
\[
 \delta^x_j=v_{X_j}(x_j),\qquad b^x_j=s(X_j),
 \qquad
 \delta^c_i=v_{S_i}((Ax)_i),\qquad b^c_i=s(S_i).
\]
Thus, $\delta^x$ measures violations of the variable bounds $x\in X$, while
$\delta^c$ measures violations of the row bounds $Ax\in S$.

\paragraph{Nine acceptance criteria.}
The checker evaluates
\begin{equation}
\begin{aligned}
 g_1&=\frac{\lVert\delta^x\rVert_\infty}
            {1+\lVert b^x\rVert_\infty}, &
 g_2&=\max_j\frac{\delta^x_j}{1+b^x_j}, &
 g_3&=\frac{\lVert\delta^c\rVert_2}
            {1+\lVert b^c\rVert_2},\\
 g_4&=\max_i\frac{\delta^c_i}{1+b^c_i}, &
 g_5&=\frac{\lVert e\rVert_2}{1+\lVert c\rVert_2}, &
 g_6&=\max_j\frac{|e_j|}{1+|c_j|},\\
 g_7&=\max_i\operatorname{dist}(y_i,D(S_i)), &
 g_8&=\max_j\operatorname{dist}(r_j,D(X_j)), &
 g_9&=\frac{|p-d|}{1+|p|+|d|}.
\end{aligned}
\label{eq:uniform-check}
\end{equation}
The pairs $(g_1,g_2)$, $(g_3,g_4)$, and $(g_5,g_6)$ measure variable-bound feasibility, row feasibility, and stationarity, respectively.  In each pair, the first quantity is a normwise residual and the second is the worst coordinatewise normalized residual.  Criteria $g_7$ and $g_8$ measure violations of the required dual and reduced-cost signs, and $g_9$ is the relative primal--dual objective gap.  The added $1$ in each denominator keeps the normalization well defined when the corresponding data or objective is zero.

A candidate is accepted if and only if every computed quantity is finite and
\[
 \max_{1\le q\le 9} g_q\le 10^{-6}.
\]
Table~\ref{tab:scale-residuals} displays $g_3$, $g_5$, and $g_9$ in the
``rel. primal,'' ``rel. stationarity,'' and ``rel. gap'' columns.  Its
``max. criterion'' column is the maximum over all nine tests, including the
six not printed separately.  The separate checker evaluates the tests in
original coordinates using long-double accumulation.

\begin{table}[!ht]
\centering
\caption{Iterations, runtimes, and separate original-space checks for the scale experiments. Convergence is evaluated every 200 iterations, giving iteration counts in multiples of 200. Times are solver/end-to-end. The maximum criterion is the largest of the nine quantities in Equation~\eqref{eq:uniform-check}.}
\label{tab:scale-residuals}
\scriptsize
\setlength{\tabcolsep}{2.4pt}
\begin{tabular}{@{}lrrrrrr@{}}
\toprule
Instance & iter. & \shortstack{solver /\\end-to-end (s)} & \shortstack{rel.\\primal} &
\shortstack{rel.\\stationarity} & \shortstack{rel.\\gap} & \shortstack{max.\\criterion} \\
\midrule
\texttt{tai256c} AJ & 3,200 & 717 / 920 & $3.62\times10^{-10}$ & $8.87\times10^{-17}$ & $5.82\times10^{-10}$ & $1.26\times10^{-9}$ \\
KDD12 L1-SVM & 74,600 & 37,294 / 37,587 & $1.82\times10^{-10}$ & $1.23\times10^{-10}$ & $4.57\times10^{-10}$ & $6.72\times10^{-7}$ \\
MS1 matching & 31,800 & 3,721 / 3,760 & $4.10\times10^{-11}$ & $1.74\times10^{-16}$ & $5.27\times10^{-15}$ & $9.10\times10^{-8}$ \\
MCF4.15B & 1,200 & 490 / 812 & $1.01\times10^{-9}$ & $3.43\times10^{-16}$ & $7.25\times10^{-9}$ & $2.86\times10^{-8}$ \\
MCF5B & 1,200 & 560 / 634 & $1.01\times10^{-9}$ & $8.21\times10^{-15}$ & $7.26\times10^{-9}$ & $2.86\times10^{-8}$ \\
MCF8B & 1,200 & 14,905 / 15,510 & $1.01\times10^{-9}$ & $2.92\times10^{-16}$ & $7.27\times10^{-9}$ & $2.86\times10^{-8}$ \\
MCF8.35B & 4,800 & 2,933 / 3,045
& $3.15\times10^{-10}$ & $1.25\times10^{-13}$
& $2.19\times10^{-9}$ & $2.19\times10^{-9}$ \\
MCF13.60B & 4,600 & 1,905 / 2,057
& $1.04\times10^{-10}$ & $1.32\times10^{-13}$
& $6.74\times10^{-8}$ & $6.74\times10^{-8}$ \\
\bottomrule
\end{tabular}
\end{table}

\arxivchange{The distributed Google benchmark runs evaluate convergence every 200 main iterations; Heat Source Easy uses an interval of 64. In the three single-GPU cases---TSP-Gaia-10M, Production Inventory, and Supply Chain---the solver uses absolute primal and stationarity tolerances of $10^{-8}$ and a relative-gap tolerance of $10^{-2}$. In these three cases, polishing is attempted after 100, 200, 400, $\ldots$ main iterations when the relative gap is at most $10^{-2}$.}

\arxivreplaced{Solver time is measured inside the iterative solver and includes residual evaluations, restart, and feasibility polishing when enabled. End-to-end time adds rank-local input, initialization, solution recovery, and output, but excludes offline conversion and the separate final check.}{Solver time \arxivgreen{in Tables~\ref{tab:google} and~\ref{tab:scale}} is the elapsed optimization time at termination. For runs evaluated under the Google study's criterion, the solver stops when a scheduled convergence evaluation satisfies that criterion or when a polished candidate satisfies the same test. This clock includes residual evaluations, restarts, unsuccessful polishing attempts, and the successful final attempt when used. End-to-end time also includes loading the prepared blocks, initialization, recovery, and output. Constructing those blocks and the subsequent independent check are excluded. The independent check determines whether the final exported solution is counted as solved.}

 In Table~\ref{tab:google}, TSP-Gaia-10M, Production Inventory, QAP-WIL-100, and Supply Chain are medians of three runs; the remaining solved rows use one checked run. Table~\ref{tab:scale} uses one checked run each for KDD12 L1-SVM, \texttt{tai256c} AJ, and MS1 matching, \arxivgreenreplace{the componentwise median of three MCF5B runs}{the median of each reported quantity across three MCF5B runs}, and the first of two checked MCF8B runs, whose end-to-end times were 15,509.9 and 15,519.7 seconds. Table~\ref{tab:scale-residuals} additionally includes one checked MCF4.15B run; it used two shared compute nodes, so its timing is descriptive rather than a scaling measurement. MCF8.35B and MCF13.60B each report one cold-start execution without presolve or feasibility polishing. Internal tolerances are $10^{-8}$ and $10^{-7}$, respectively; both pass the same independent $10^{-6}$ check.

\arxivchange{Heat Source Easy uses a $1\times8$ grid on eight Quadro RTX~6000 nodes, without presolve, with internal feasibility and relative-gap tolerances of $10^{-6}$. It terminates after 935,040 iterations and 28,858.4 seconds of solver time; the end-to-end time is 28,877.9 seconds, excluding the separate final check. All nine original-space criteria are at most $9.26\times10^{-7}$.}

The Google study used a six-day limit \citep{applegate2026pdlp}; accordingly, Table~\ref{tab:google} compares coverage and reports raw runtimes rather than a controlled hardware comparison.

The published CPU experiments used either a 16-core AMD EPYC 7302 compute node with 256 GB of memory or a 32-core Intel Xeon Platinum 8352Y compute node with 1 TB. PDLP used 16 or 32 threads, while Gurobi 11.0.2 used its default thread count and termination rule.  Table~\ref{tab:google} reports the fastest published Gurobi time among barrier, primal simplex, and dual simplex \citep{applegate2026pdlp}.

MS1 uses a $1\times32$ layout over sixteen nodes with two V100 GPUs each (16- and 32-GiB variants), socket transport, a cold start, and no presolve or feasibility polishing. The exported solution passes both the original-LP checker and an independent original-source audit at $10^{-6}$. Its end-to-end time is the solver-process wall time, excluding the separate final check. MCF8B uses a $1\times 12$ layout over three compute nodes, without presolve or feasibility polishing. Hardware, software versions, interconnect settings, and allocation details for all experiments are reported in Appendix~\ref{app:hardware}.

\FloatBarrier

\arxivchange{\section{Memory layout and workspace reuse}
\label{app:memory}}

\arxivchange{\subsection{Sparse matrix and index storage}}

\arxivchange{The local matrix must support both $A_{rc}x_c$ and $A_{rc}^{\top}y_r$. In compressed sparse row (CSR) form, the entries of a row are stored consecutively, with a pointer to each row's start. We store CSR forms of both $A_{rc}$ and $A_{rc}^{\top}$, providing this row-wise layout for both products. Let $m_r=|I_r|$, $n_c=|J_c|$, and $z_{rc}=\nnz(A_{rc})$. Each orientation stores an 8-byte value and a 4-byte index per nonzero. Together with the two row-pointer arrays, these matrix arrays occupy}
\arxivchange{\begin{equation}
 M_A=24z_{rc}+4(m_r+n_c+2)\quad\text{bytes}.
 \label{eq:matrix-storage}
\end{equation}}
\arxivchange{This is the storage for the two sparse matrices; vectors and library workspaces require additional memory. Global coordinates are recovered by adding a block's 64-bit offset to its local index. Thus, the complete LP can exceed 32-bit counts while each local row, column, and nonzero count must remain within the signed 32-bit range.}

\arxivgreen{The last CSR row pointer equals the stored nonzero count. Keeping these pointers local avoids imposing a signed 32-bit limit on the complete matrix.}

\arxivchange{\subsection{Bounds and temporary vectors}}

\arxivchange{For a block of nonnegative variables with no finite upper bounds, projecting a vector $v$ reduces to $v_j\mapsto\max\{v_j,0\}$. Storing a zero lower bound and an infinite upper bound for every variable is unnecessary. The implementation checks the scaled bounds and omits each constant array independently; other bounds retain their arrays. Each omitted FP64 array saves $8n_c$ bytes. With 100 million local variables, omitting both arrays saves 1.6 decimal GB without changing the feasible set or the projection.}

\arxivchange{The same storage reduction is possible when temporary values are used at different times. On the primal side, a reflected iterate is first read by the matrix product. Its storage is then available for the iterate difference used in the fixed-point calculation, and subsequently for the stationarity residual. Once the last reader of one value has finished, the next value can overwrite it. One allocation therefore replaces three, saving $16n_c$ bytes compared with storing these three vectors separately. The row-product buffer is reused similarly for dual differences and row residuals. Their evaluation order ensures that a value is consumed before its storage is overwritten.}

\arxivchange{These savings depend on the local slice lengths, not just on the number of nonzeros. Primal vectors have $n_c$ entries and are replicated down a process column; dual vectors have $m_r$ entries and are replicated across a process row. A grid chosen only to balance matrix blocks can therefore leave too little space for the vectors. Both contributions must be included when choosing the partition.}

\arxivchange{\subsection{Memory used by feasibility polishing}}

\arxivchange{Feasibility polishing solves separate primal and dual feasibility subproblems. \arxivgreen{The matrix block is shared by the main solve and both phases.} Each requires additional iterates and work vectors alongside the main PDHG state. Keeping the working vectors for both subproblems on the GPU can raise peak memory even when the main iteration fits.}

\arxivchange{Their work vectors can instead be allocated one phase at a time. After one phase succeeds, its local candidate is copied to host memory and its GPU work vectors are freed before the other phase begins. A dual candidate retains the reduced-cost vector computed with it. The main iterate is preserved while polishing is attempted, so a failed attempt can resume the original solve. This exchanges host storage and transfers for fewer simultaneous GPU allocations. \arxivgreen{Figure~\ref{fig:polishing-memory} compares the two allocation patterns.}}

\arxivgreen{\begin{figure}[!ht]
  \centering
  {\color{black}\resizebox{\linewidth}{!}{
\begingroup
\definecolor{memblue}{RGB}{0,114,178}
\definecolor{memorng}{RGB}{210,140,0}
\definecolor{memteal}{RGB}{0,133,119}
\begin{tikzpicture}[
  x=1cm,y=1cm,font=\fontsize{9}{10.5}\selectfont,
  panel/.style={draw=black!23,fill=black!1,rounded corners=2pt,line width=.45pt},
  flow/.style={-{Latex[length=1.8mm,width=1.1mm]},draw=memblue,line width=.65pt},
  labeltext/.style={align=center,inner sep=1pt},
  retained/.style={draw=black!35,fill=black!7,line width=.45pt}
]
\path[use as bounding box] (0,-.62) rectangle (15.8,7.3);
\draw[panel] (0,0) rectangle (7.76,7.3);
\draw[panel] (8.04,0) rectangle (15.8,7.3);
\node[anchor=west,font=\fontsize{10.3}{12}\selectfont\bfseries] at (.20,6.91)
  {(a) Retain both phase workspaces};
\node[anchor=west,font=\fontsize{10.3}{12}\selectfont\bfseries] at (8.24,6.91)
  {(b) Allocate one phase at a time};

\foreach \base in {0,8.04} {
  \begin{scope}[shift={(\base,0)}]
  \node[labeltext] at (3.37,6.24) {Primal phase};
  \node[labeltext] at (6.00,6.24) {Dual phase};
  \draw[-{Latex[length=1.5mm]},draw=black!60] (2.12,5.94) -- (7.28,5.94);
  \draw[draw=black!18,fill=white,rounded corners=2pt]
    (.18,2.10) rectangle (7.58,5.70);
  \node[anchor=west,font=\scriptsize\bfseries,text=black!65] at (.34,5.43) {GPU};
  \node[labeltext,anchor=west,text width=1.54cm] at (.34,4.95) {Matrix};
  \node[labeltext,anchor=west,text width=1.54cm] at (.34,4.18) {Main\\state};
  \node[labeltext,anchor=west,text width=1.54cm] at (.34,3.38) {Primal\\workspace};
  \node[labeltext,anchor=west,text width=1.54cm] at (.34,2.58) {Dual\\workspace};
  \draw[retained] (2.12,4.68) rectangle (7.28,5.22);
  \node at (4.70,4.95) {Shared matrix block};
  \draw[draw=memorng!70!black,fill=memorng!14,line width=.45pt]
    (2.12,3.91) rectangle (7.28,4.45);
  \node at (4.70,4.18) {Preserved for resuming the main solve};
  \draw[draw=memblue,fill=memblue!16,line width=.45pt]
    (2.12,3.11) rectangle (4.62,3.65);
  \node[text=memblue!85!black] at (3.37,3.38) {active};
  \draw[draw=memteal,fill=memteal!17,line width=.45pt]
    (4.75,2.31) rectangle (7.28,2.85);
  \node[text=memteal!90!black] at (6.00,2.58) {active};
  \end{scope}
}

\draw[draw=memblue,fill=memblue!7,line width=.45pt]
  (4.62,3.11) rectangle (7.28,3.65);
\node[text=memblue!85!black] at (6.00,3.38) {still allocated};
\node[labeltext,text width=6.4cm] at (3.88,1.30)
  {During dual polishing, both sets of\\phase work vectors occupy GPU memory.};
\node[font=\bfseries,text=black!75] at (3.88,.45)
  {Two phase workspaces coexist};

\draw[draw=memblue!45,dashed,fill=white,line width=.45pt]
  (12.79,3.11) rectangle (15.32,3.65);
\node[text=black!50] at (14.04,3.38) {freed};
\draw[flow] (12.58,3.10) -- (12.58,1.49);
\node[anchor=east,font=\scriptsize,text=memblue!85!black] at (12.45,1.77)
  {copy, then free};
\draw[draw=memblue!65,fill=memblue!8,rounded corners=2pt,line width=.45pt]
  (10.16,.83) rectangle (15.32,1.49);
\node[labeltext,text=memblue!85!black] at (12.74,1.16)
  {Primal candidate retained on host};
\node[font=\bfseries,text=memteal!90!black] at (11.92,.45)
  {One GPU phase workspace at a time};
\node[font=\fontsize{8}{9}\selectfont,text=black!65] at (7.9,-.34)
  {Reported distributed Google runs use (a); (b) illustrates sequential allocation.};
\end{tikzpicture}
\endgroup}}
  \caption{\arxivgreen{GPU storage during feasibility polishing. (a) Retaining both phase workspaces raises memory use during the second phase. (b) Copying a successful candidate to host memory lets the first workspace be freed before the next is allocated. The matrix is shared in both cases. The reported distributed Google runs use (a).}}
  \label{fig:polishing-memory}
\end{figure}}

\arxivchange{The distributed Google benchmark runs use the bound and workspace reductions above, but their polishing implementation retains both sets of working vectors on the GPU. Their timings therefore do not measure the effect of sequential allocation. Heat Source Easy uses workspace reuse together with iterate averaging and no feasibility polishing. Table~\ref{tab:storage-memory} reports sampled device-memory peaks for the recorded runs; those measurements are not separate comparisons of the individual storage choices.}

\arxivgreen{\FloatBarrier}

\section{Support-aware update equivalence and communication model}
\label{app:participant}

For a $1\times p$ column partition, define the participant set of row $i$ by
\begin{equation}
\mathcal P_i :=
\left\{j\in\{0,\ldots,p-1\} :
\nnz(A_{i,J_j})>0\right\},
\qquad
k_i := |\mathcal P_i|.
\label{eq:participants}
\end{equation}
Thus, $\mathcal P_i$ contains exactly the GPUs whose column shards store a nonzero coefficient in row $i$.

\begin{figure}[!ht]
  \centering
  \resizebox{.98\linewidth}{!}{\begingroup
\definecolor{rankblue}{RGB}{0,114,178}
\definecolor{rankorange}{RGB}{230,159,0}
\definecolor{accentteal}{RGB}{0,133,119}
\begin{tikzpicture}[
  x=1cm,y=1cm,
  font=\footnotesize,
  panel/.style={draw=black!24,fill=black!1,rounded corners=2pt,line width=.45pt},
  cell/.style={draw=black!45,rounded corners=1.2pt,minimum width=.86cm,
               minimum height=.46cm,inner sep=1.2pt,align=center},
  op/.style={rounded corners=1.5pt,minimum width=2.75cm,
             minimum height=.60cm,inner sep=2pt,align=center},
  flow/.style={-{Latex[length=1.8mm,width=1.1mm]},line width=.5pt,draw=black!70}
]

\draw[panel] (0,0) rectangle (4.82,3.65);
\node[anchor=west,font=\small\bfseries] at (.22,3.35) {(a) Row incidence};
\node[anchor=west] at (.35,2.92) {column shard};
\foreach \x/\j/\fillc in {1.18/0/rankblue!12,2.18/1/black!3,3.18/2/rankorange!16,4.18/3/black!3} {
  \node[cell,fill=\fillc] at (\x,2.45) {$j=\j$};
}
\node[cell,fill=rankblue!12]   (a0) at (1.18,1.72) {$A_{i,J_0}$};
\node[cell,fill=black!3,text=black!45] (a1) at (2.18,1.72) {$0$};
\node[cell,fill=rankorange!16] (a2) at (3.18,1.72) {$A_{i,J_2}$};
\node[cell,fill=black!3,text=black!45] (a3) at (4.18,1.72) {$0$};
\node[cell,fill=rankblue!12]   at (1.18,.90) {$a_i^{(0)}$};
\node[cell,fill=black!3,text=black!45] at (2.18,.90) {$0$};
\node[cell,fill=rankorange!16] at (3.18,.90) {$a_i^{(2)}$};
\node[cell,fill=black!3,text=black!45] at (4.18,.90) {$0$};
\foreach \x in {1.18,2.18,3.18,4.18}
  \draw[flow] (\x,1.47) -- (\x,1.17);
\node[font=\bfseries,text=accentteal] at (2.68,.28)
  {$\mathcal P_i=\{0,2\}$};

\draw[panel] (5.08,0) rectangle (10.20,3.65);
\node[anchor=west,font=\small\bfseries] at (5.30,3.35)
  {(b) Dense update};
\foreach \x/\lab/\fillc in {5.68/{a_i^{(0)}}/rankblue!12,6.78/0/black!3,
                            8.50/{a_i^{(2)}}/rankorange!16,9.60/0/black!3} {
  \node[cell,fill=\fillc] (d\x) at (\x,2.62) {$\lab$};
}
\node[op,draw=rankblue,fill=rankblue!7] (dense) at (7.64,1.70)
  {$s_i=\sum_{j=0}^{3}a_i^{(j)}$\\each rank computes $y_i^+=\Phi_i(y_i,s_i)$};
\foreach \x in {5.68,6.78,8.50,9.60}
  \draw[flow] (\x,2.38) -- (dense.north);
\foreach \x/\fillc in {5.68/rankblue!12,6.78/black!3,
                       8.50/rankorange!16,9.60/black!3} {
  \node[cell,fill=\fillc] (dy\x) at (\x,.66) {$y_i^+$};
  \draw[flow] (dense.south) -- (\x,.91);
}
\node[font=\scriptsize,text=black!55] at (7.64,.18)
  {communicate and store on all four ranks};

\draw[panel] (10.46,0) rectangle (15.58,3.65);
\node[anchor=west,font=\small\bfseries] at (10.68,3.35)
  {(c) Participant update};
\node[cell,fill=rankblue!12]   (p0) at (11.18,2.62) {$a_i^{(0)}$};
\node[cell,fill=black!2,draw=black!18,text=black!35] at (12.28,2.62) {--};
\node[cell,fill=rankorange!16] (p2) at (14.00,2.62) {$a_i^{(2)}$};
\node[cell,fill=black!2,draw=black!18,text=black!35] at (15.10,2.62) {--};
\node[op,draw=accentteal,fill=accentteal!7] (part) at (13.08,1.70)
  {owner: $s_i=a_i^{(0)}+a_i^{(2)}$\\$y_i^+=\Phi_i(y_i,s_i)$};
\draw[flow,draw=rankblue] (p0.south) -- (part.north west);
\draw[flow,draw=rankorange] (p2.south) -- (part.north east);
\node[cell,fill=rankblue!12]   (py0) at (11.18,.66) {$y_i^+$};
\node[cell,fill=black!1,draw=black!18,text=black!35] at (12.28,.66) {--};
\node[cell,fill=rankorange!16] (py2) at (14.00,.66) {$y_i^+$};
\node[cell,fill=black!1,draw=black!18,text=black!35] at (15.10,.66) {--};
\draw[flow,draw=rankblue] (part.south west) -- (py0.north);
\draw[flow,draw=rankorange] (part.south east) -- (py2.north);
\node[font=\scriptsize,text=black!55] at (13.08,.18)
  {communicate only on $\mathcal P_i$};

\end{tikzpicture}
\endgroup}
  \caption{Dense and support-aware updates for a row supported on column
  shards 0 and 2.  The dense path reduces the row activity and stores the
  updated dual coordinate on every rank.  The participant path reduces at
  an owner and returns the update only to ranks whose local $A^\top y$
  product uses that coordinate. \arxivgreen{Here $\Phi_i$ denotes the dual update for row $i$.}}
  \label{fig:participant-update}
\end{figure}

The reported participant runs use a $1\times p$ grid and apply the construction to the $Ax$ reduction and dual update; \arxivgreenreplace{$A^\top y$ follows its unchanged path}{the $A^\top y$ reduction is unchanged}. The communication plan is built once after all support-changing preprocessing, from scaled local compressed sparse row (CSR) support, and checked for agreement across the row communicator.
\arxivgreenreplace{Rows with the same participant set are packed together.  Full-participation buckets use NCCL Reduce/Broadcast, whereas partial buckets use}{Rows with the same participant set are packed into communication groups. If every rank participates, a group uses NCCL Reduce/Broadcast; otherwise it uses} grouped NCCL point-to-point transfers through the lowest participating rank; rank zero owns structurally empty rows.  The sparsity pattern remains fixed during a solve, so a structural change or repartitioning would require rebuilding the plan. Extending this path to a general $R\times C$ grid remains future work.

\paragraph{Exact-arithmetic equivalence.}
\arxivgreenreplace{Assume exact arithmetic and that both executions use the same arithmetic, parameters, and restart, evaluation, checkpoint, and polishing rules.}{Assume that the dense and support-aware executions use exact arithmetic, the same parameters, and the same rules for restarts, evaluations, checkpoint selection, and polishing.}  \arxivgreenreplace{At the start of each participant epoch, the dense and participant executions must have the same fully materialized solver state---including current, average, anchor, and reflected primal--dual components.}{At a full-state synchronization, suppose the two executions have identical current, average, anchor, and reflected primal--dual states.}  The support of $A$ must then remain fixed, every $\mathcal P_i$ in Equation~\eqref{eq:participants} must be the exact row support over ranks, and each nonempty set must have an owner $o_i\in\mathcal P_i$ holding the authoritative $y_i$.

For $j\notin\mathcal P_i$, $A_{i,J_j}=0$.  Let $\widetilde x$ denote the vector supplied to the $Ax$ product (for Equation~\eqref{eq:pdhg-dual}, $\widetilde x=2\widehat x-x$), and define $a_i^{(j)}=A_{i,J_j}\widetilde x_{J_j}$; then $a_i^{(j)}=0$ whenever $j\notin\mathcal P_i$.  Thus $\sum_j a_i^{(j)}=\sum_{j\in\mathcal P_i}a_i^{(j)}$, which gives the same scalar row activity to the $i$th coordinate of $\Phi_{S,\sigma}$.  The owner computes the same coordinate update, and dissemination reaches every rank whose transpose product contains that coordinate.  The local Halpern and vector updates therefore match.

Before evaluation, adaptive restart, checkpoint selection, polishing, or any other operation requiring replicated dual state, the participant execution reconstructs every required coordinate; global statistics count each authoritative coordinate once.  \arxivgreenreplace{Induction over these epochs then proves whole-solver equality at synchronization points.}{Starting from identical initial states, induction over the synchronization intervals gives agreement at each full-state synchronization.}  If $\mathcal P_i=\varnothing$, activity is identically zero and a designated owner evaluates the $i$th dual update with zero row activity and no operator communication.

Floating-point reductions can use different summation orders.  Across the shared evaluation checkpoints of the V100 Design Match paired runs, the largest absolute differences between the reported dense and participant fields are $1.87\times10^{-12}$ for relative primal residual, $1.98\times10^{-11}$ for relative stationarity, and $1.30\times10^{-15}$ for relative gap.

For $p>1$, count one scalar transferred across one rank-to-rank hop as one logical scalar-hop.  Let $\mathcal S$ be the $m_s$ rows assigned to the participant path, $h$ the number of dual updates, $b$ the number of full-state synchronization boundaries, and $q$ the number of dual vectors synchronized at each boundary.  Under the ring model, the corresponding counts are
\begin{equation}
 H_{\mathrm{dense}}=2h(p-1)m_s,\qquad
 H_{\mathrm{part}}=2h\sum_{i\in\mathcal S}(k_i-1)_+ + qb(p-1)m_s,
 \qquad
\rho_{\mathrm{tot}}=1-\frac{H_{\mathrm{part}}}{H_{\mathrm{dense}}},
\label{eq:traffic-reduction}
\end{equation}
where $(t)_+=\max\{t,0\}$.  The factor of two counts the reduction and dissemination legs of each dual update; the boundary term counts a one-way full-state synchronization.  Figure~\ref{fig:med-participant}b reports reductions in this logical scalar-hop count, not measured bytes on physical links.  Packing, imbalance, topology, and latency are excluded, so the model need not predict elapsed time.  For Design Match, $p=8$, $m_s=22{,}000{,}135$, $h=27{,}081$, $b=138$, $q=1$, and $\sum_i(k_i-1)_+=10{,}426{,}445$, giving the reported 92.9748\%.  Hypergraph partitioning, configured sparse collectives, and node-aware sparse matrix--vector multiplication provide the related sparse-communication work \citep{ucar2007hypergraph,zhao2013sparse,hoefler2009sparse,hough2026sparse,bienz2019node}.  \arxivgreenreplace{We do not claim a new generic collective: the contribution is the solver-specific integration of a fixed producer/consumer plan with the nonlinear $\Phi_{S,\sigma}$ update, empty-row ownership, reflected and restart state, polishing, checking, and output boundaries.}{\solvername{} combines these sparse exchanges with the dual update, restoring the replicated state needed for convergence checks, restarts, polishing, and output.}


\section{Related work and comparison scope}
\label{app:scope}

For this comparison, we focus on GPU LP solvers that use a general sharded sparse-matrix representation and solver kernels, rather than application-specific operators or decompositions. \arxivgreenreplace{We count a run only if it exports a complete primal--dual solution that passes the independent checker.}{For \solvername{}, we count only runs that export a complete primal--dual solution and pass the independent checker.}

D-PDLP evaluates its two-dimensional decomposition on a single compute node containing eight H100 GPUs connected by NVLink.  Its largest-variable instance has 127.5 million variables and 257.5 million nonzeros; another reaches 690 million nonzeros \citep{li2026dpdlp}.  MPAX reports dense multi-GPU PDHG at approximately 900 million represented nonzeros on up to four H100 GPUs, with efficient distributed sparse-data sharding left as future work \citep{lu2024mpax}.  NVIDIA's cuOpt 26.08 release notes report $2.5\times$--$8.8\times$ PDLP speedups on eight NVLink-connected B200 GPUs; they do not identify a cross-node execution.\footnote{\href{https://docs.nvidia.com/cuopt/user-guide/latest/release-notes.html}{NVIDIA cuOpt 26.08 release notes}.} The cuOpt FAQ gives single-H100 capacity examples up to two billion nonzeros.\footnote{\href{https://docs.nvidia.com/cuopt/user-guide/latest/faq.html}{NVIDIA cuOpt FAQ}.}

A specialized regularized-matching solver spans sixteen H100 GPUs across two compute nodes and reports models up to one billion nonzeros \citep{rahmattalabi2026matching}.  It uses a source-decomposable formulation and a dual-gradient method specialized to matching.  ECLIPSE reports a $10^{12}$-variable structured web LP by optimizing the smooth dual of a ridge-perturbed model with distributed matrix--vector products \citep{basu2020eclipse}. 

Among the general sparse GPU LP solvers discussed above, MCF13.60B is, to our knowledge, the largest reported validated solve, with 40.807 billion nonzeros. This is a scale claim, not a claim of being the first multi-node GPU optimization method or of a matched speedup over the systems above. Additionally, the 70-GPU capacity experiment fits a deterministic flow model with 8.500 billion variables, 107.374 million constraints, and 17.045 billion nonzeros across twenty compute nodes.  It ran 280 iterations in 217.46 seconds with 339~MiB free on the fullest GPU, but did not export the complete vectors needed for a separate primal--dual check.  The time characterizes this capacity execution, not a checked solve or a time-to-solution result.




\end{document}